\documentclass[12pt]{article}
\usepackage{graphicx}
\usepackage{pstcol}
\font\bg=cmbx10 scaled\magstep1
\font\Bg=cmbx12 scaled\magstep3
\font\small=cmr8

\newtheorem{newlemma}{{\bf Lemma}}

\newtheorem{newteorem}{{\bf Theorem}}

\newenvironment{teorem}{\begin{newteorem}{\hspace{-0.5
em}{\bf.}}}{\end{newteorem}}

\newtheorem{newkorolari}{{\bf Corollary}}

\newtheorem{newdefine}{{\bf Definition}}

\newtheorem{newquestion}{{\bf Question}}

\newtheorem{newkonjek}{{\bf Conjecture}}

\newtheorem{newexample}{{\bf Example}}

\begin{document}
\tolerance=10000
\baselineskip18truept
\newbox\thebox
\global\setbox\thebox=\vbox to 0.2truecm{\hsize
0.15truecm\noindent\hfill}
\def\boxit#1{\vbox{\hrule\hbox{\vrule\kern0pt
     \vbox{\kern0pt#1\kern0pt}\kern0pt\vrule}\hrule}}
\def\qed{\lower0.1cm\hbox{\noindent \boxit{\copy\thebox}}\bigskip}
\def\ss{\smallskip}
\def\ms{\medskip}
\def\bs{\bigskip}
\def\c{\centerline}
\def\nt{\noindent}
\def\ul{\underline}
\def\ol{\overline}
\def\lc{\lceil}
\def\rc{\rceil}
\def\lf{\lfloor}
\def\rf{\rfloor}
\def\ov{\over}
\def\t{\tau}
\def\th{\theta}
\def\k{\kappa}
\def\l{\lambda}
\def\L{\Lambda}
\def\g{\gamma}
\def\d{\delta}
\def\D{\Delta}
\def\e{\epsilon}
\def\lg{\langle}
\def\rg{\rangle}
\def\p{\prime}
\def\sg{\sigma}
\def\ch{\choose}

\newcommand{\ben}{\begin{enumerate}}
\newcommand{\een}{\end{enumerate}}
\newcommand{\bit}{\begin{itemize}}
\newcommand{\eit}{\end{itemize}}
\newcommand{\bea}{\begin{eqnarray*}}
\newcommand{\eea}{\end{eqnarray*}}
\newcommand{\bear}{\begin{eqnarray}}
\newcommand{\eear}{\end{eqnarray}}

\centerline{\Bg  An atlas of domination polynomials}
 \vspace{.3cm}

\centerline {\Bg of graphs of order at most six}
\bigskip

\baselineskip12truept
\centerline{\bg Saeid Alikhani$^{a,}${}\footnote{\baselineskip12truept\it\small
Corresponding author. E-mail: alikhani206@gmail.com} and Yee-hock Peng$^{b,c}$}
\baselineskip20truept
\centerline{\it $^{a}$Department of Mathematics, Yazd University}
\vskip-8truept
\centerline{\it 89195-741, Yazd, Iran}
\vskip-9truept
\centerline{\it $^{b}$Institute for Mathematical Research, and}
\vskip-9truept
\centerline{\it $^{c}$Department of Mathematics,}
\vskip-8truept
\centerline{\it University Putra Malaysia, 43400 UPM Serdang, Malaysia}
\vskip-8truept
\centerline{\it }
\vskip-0.2truecm
\nt\rule{16cm}{0.1mm}

\nt{\bg ABSTRACT}
\medskip

\noindent{\it  The domination polynomial of a graph $G$ of order $n$ is the
polynomial $D(G,x)=\sum_{i=\gamma(G)}^{n} d(G,i) x^{i}$, where $d(G,i)$ is
the number of dominating sets  of $G$ of size $i$, and $\gamma(G)$ is the
domination number of $G$.  The roots of domination polynomial is called
domination roots. In this article, we compute the domination polynomial and domination roots of all graphs
of  order less than or equal to 6, and show them in the tables.}
\nt\rule{16cm}{0.1mm}

 \section{Introduction}

\nt
Let $G=(V,E)$ be a graph of order $|V|=n$. For any vertex $v\in V$,
the {\it open neighborhood} of $v$ is the set $N(v)=\{u \in V|uv\in E\}$ and
the {\it closed neighborhood} of $v$ is the set $N[v]=N(v)\cup \{v\}$.
 For a set $S\subseteq V$, the open neighborhood of $S$
 is $N(S)=\bigcup_{v\in S} N(v)$ and the closed neighborhood of $S$ is $N[S]=N(S)\cup S$.
A set $S\subseteq V$ is a dominating set if $N[S]=V$, or equivalently,
 every vertex in $V-S$ is adjacent to at least one vertex in $S$.
 The domination number $\gamma(G)$ is the minimum cardinality of a dominating set in $G$.
For a detailed treatment of this parameter, the reader is referred to~\cite{domination}.
Let ${\cal D}(G,i)$ be the family of dominating sets of a graph $G$
with cardinality $i$ and let $d(G,i)=|{\cal D}(G,i)|$.
 The {\it domination polynomial} $D(G,x)$ of $G$ is defined as
$D(G,x)=\sum_{i=\gamma(G)}^{|V(G)|} d(G,i) x^{i}$,
where $\gamma(G)$ is the domination number of $G$ (\cite{saeid1}).
 A root of $D(G,x)$ is called a {\it domination  root} of $G$.

 \ms

\nt We can obtain some information about a graph, by using its domination polynomial. (See~\cite{euro}).
Similar to another graph polynomials, it is natural that, we have a atlas for domination polynomial
and domination roots of graphs with order at most $6$. For more information on these subjects refer to \cite{euro,Oper,gcom,saeid3,few}. 

\ms

\nt In this paper, we obtain the domination polynomials and domination roots of all connected
graphs of order one to six. Let us to say that, these tables has published in \cite{thesis,book}. Thanks to ArXiv for letting us 
to publish them for more access.  Note that for disconnected graphs, we can use
the following theorem:

\begin{teorem}\label{theorem1}{\rm (\cite{saeid1})}
If a graph $G$ has  $m$ components $G_1,\ldots,G_m$, then $D(G,x)=D(G_1,x)\cdots D(G_m,x)$.
\end{teorem}

\begin{figure}[h]
\hglue.75cm
\includegraphics[width=14cm,height=10.3cm]{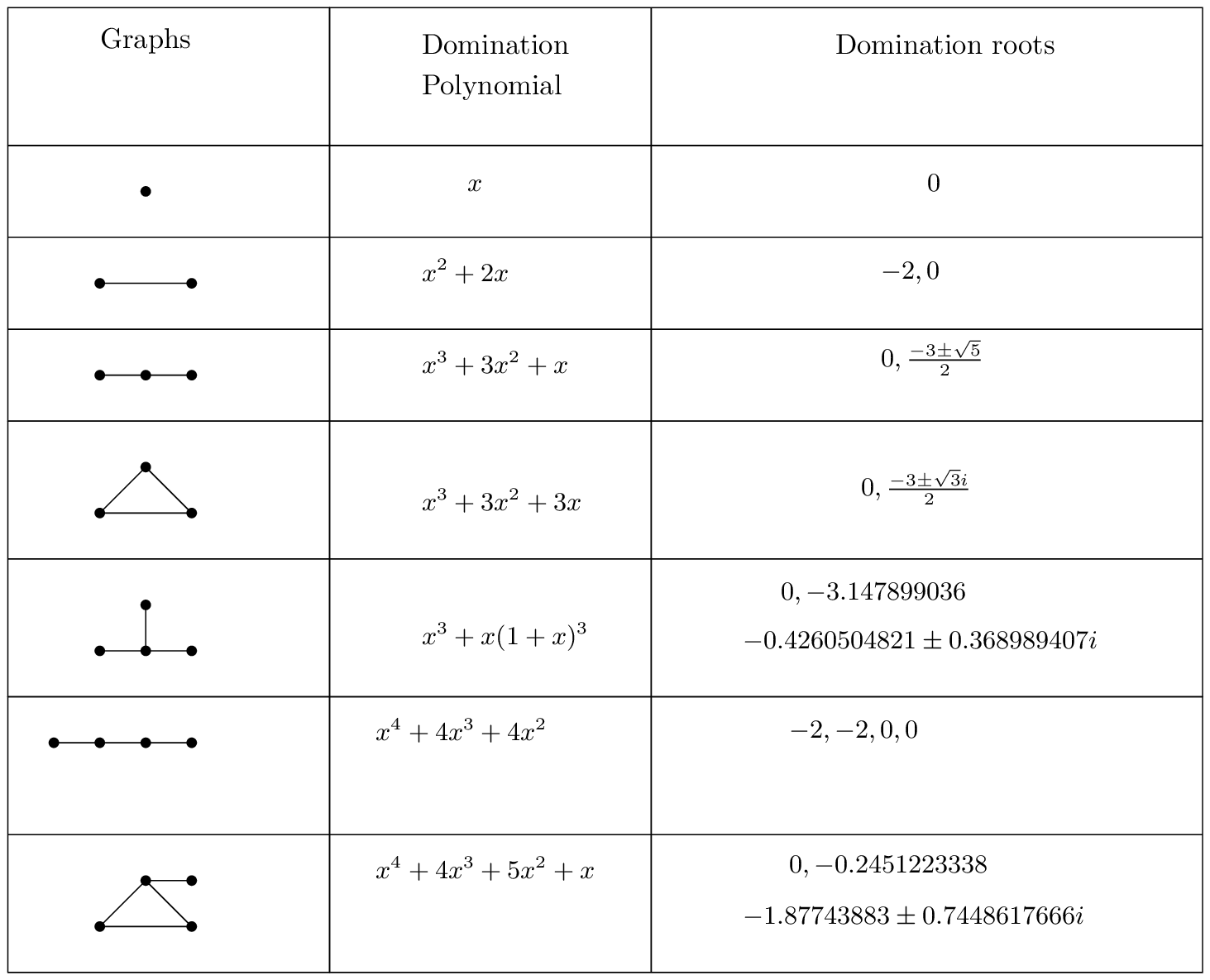}
\hglue0.75cm
\includegraphics[width=14cm,height=10.3cm]{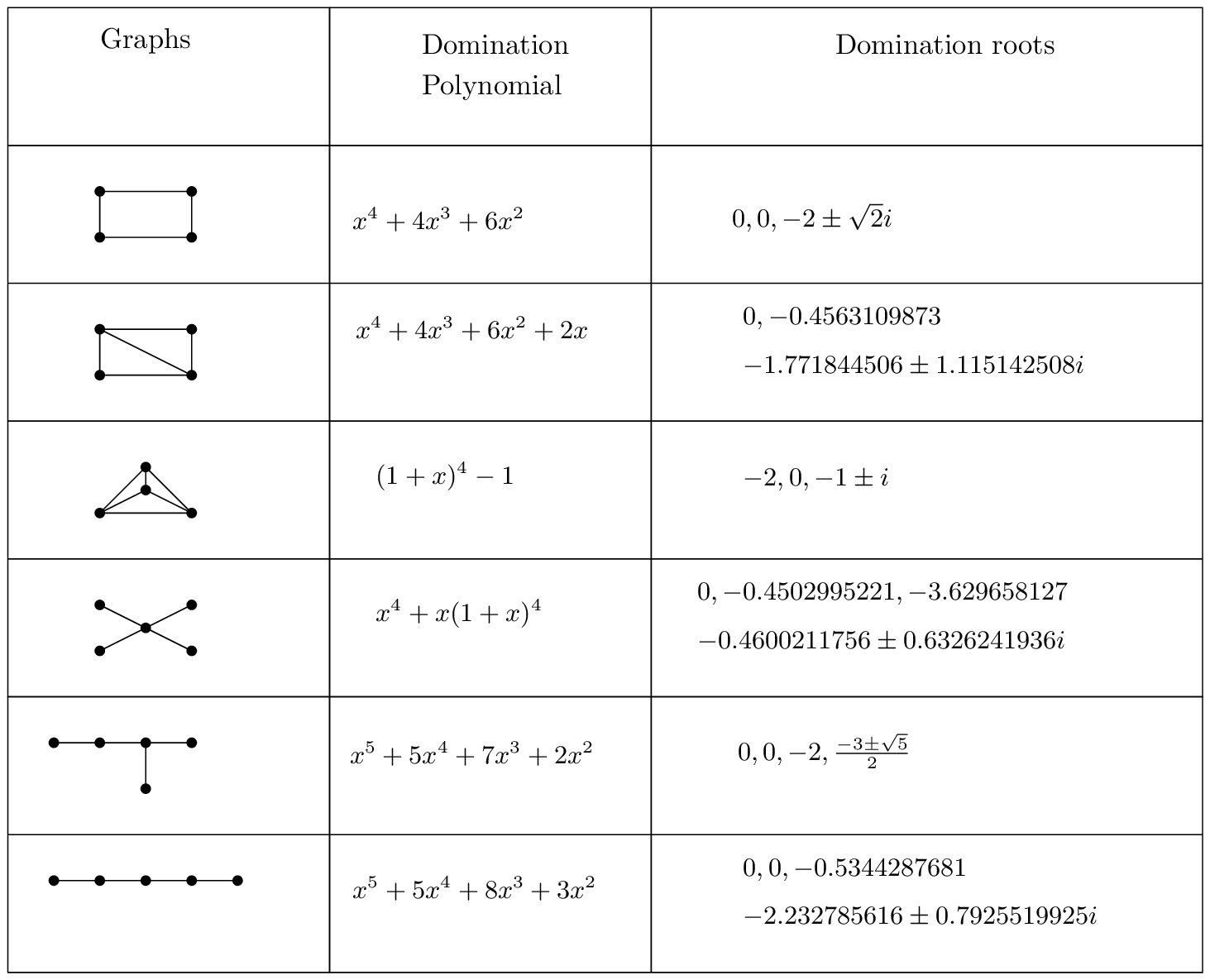}
\end{figure}

\begin{figure}[h]
\hglue0.75cm
\includegraphics[width=14cm,height=10.3cm]{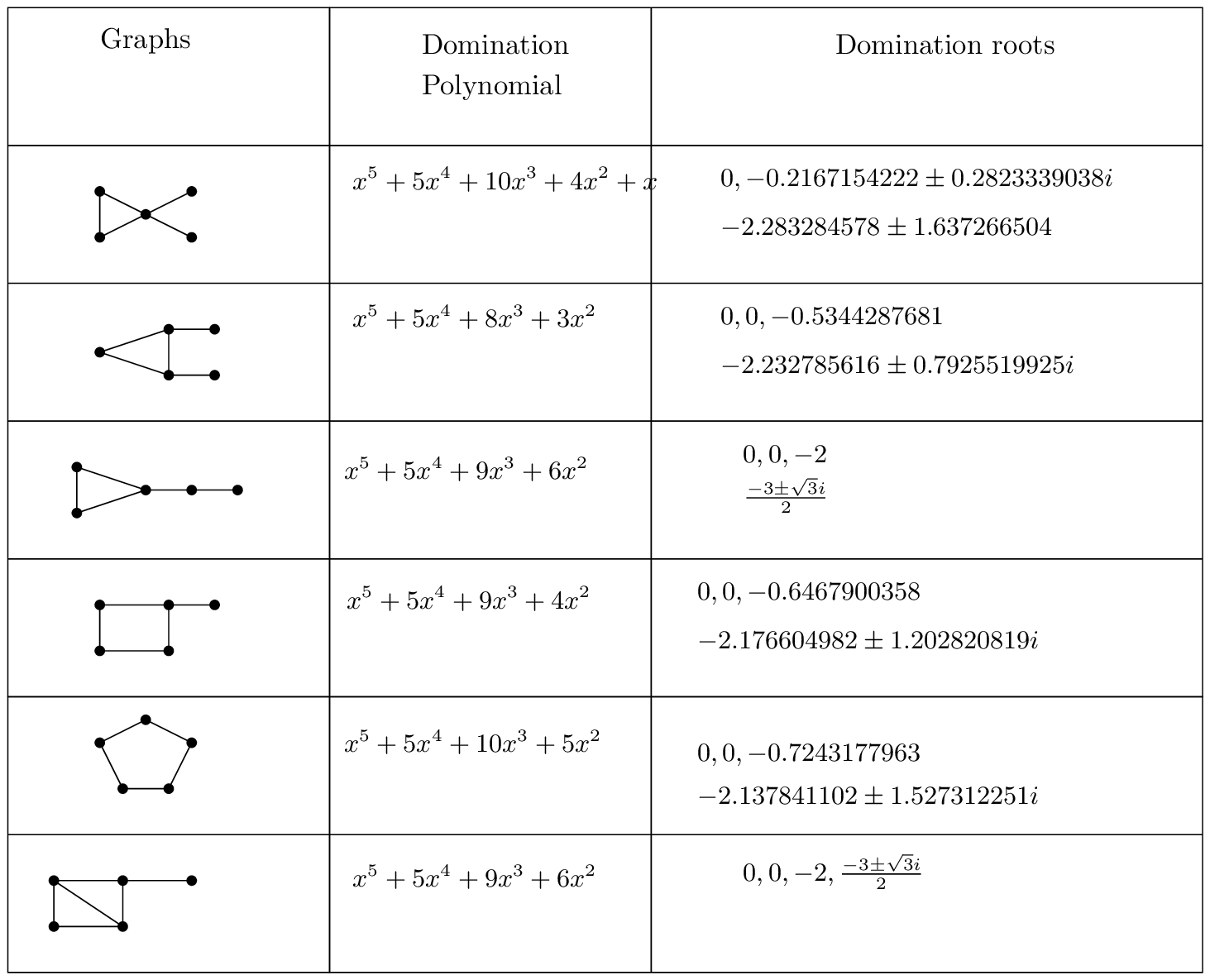}
\hglue0.75cm
\includegraphics[width=14cm,height=10.3cm]{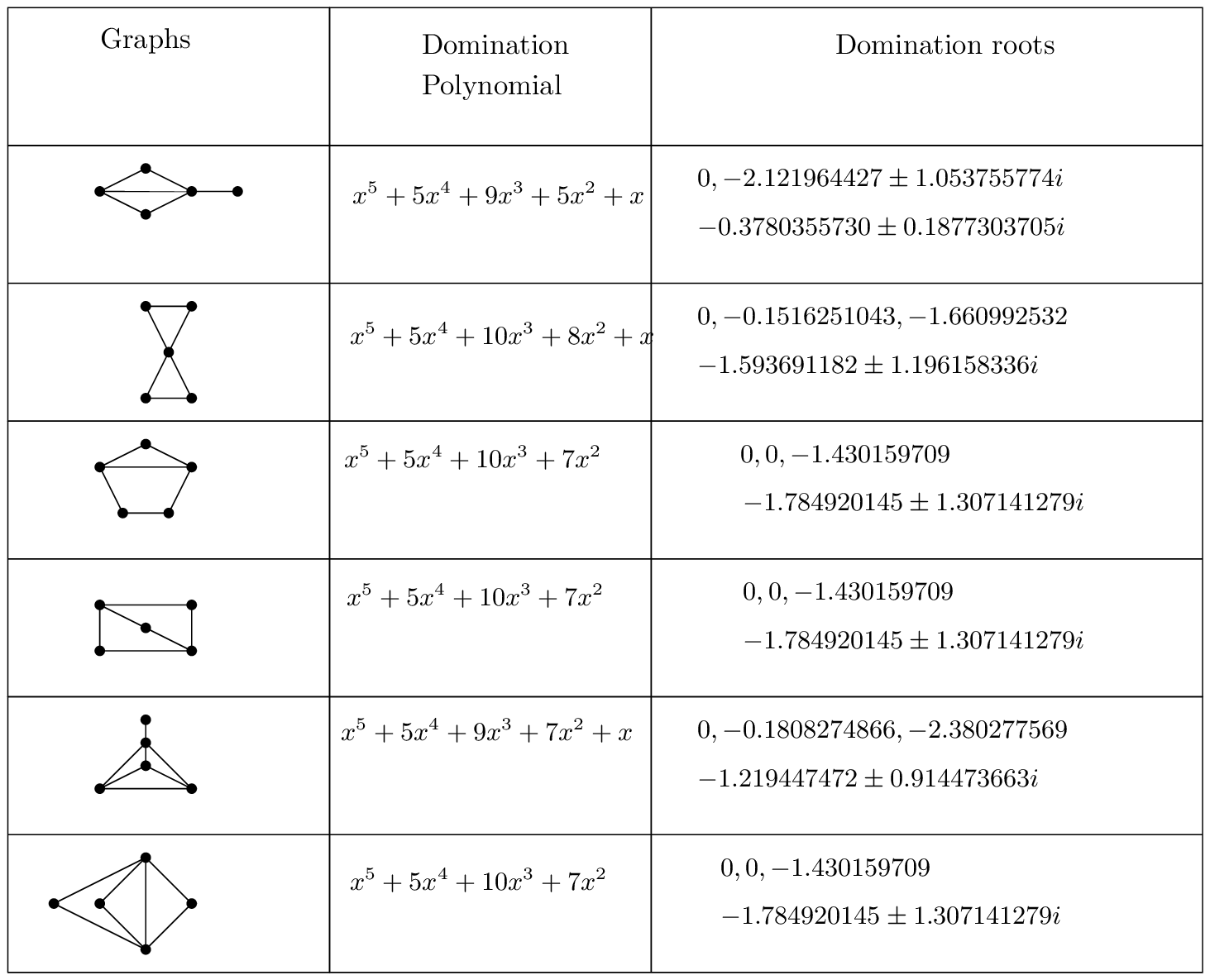}
\end{figure}

\begin{figure}[h]
\hglue0.75cm
\includegraphics[width=14cm,height=10.3cm]{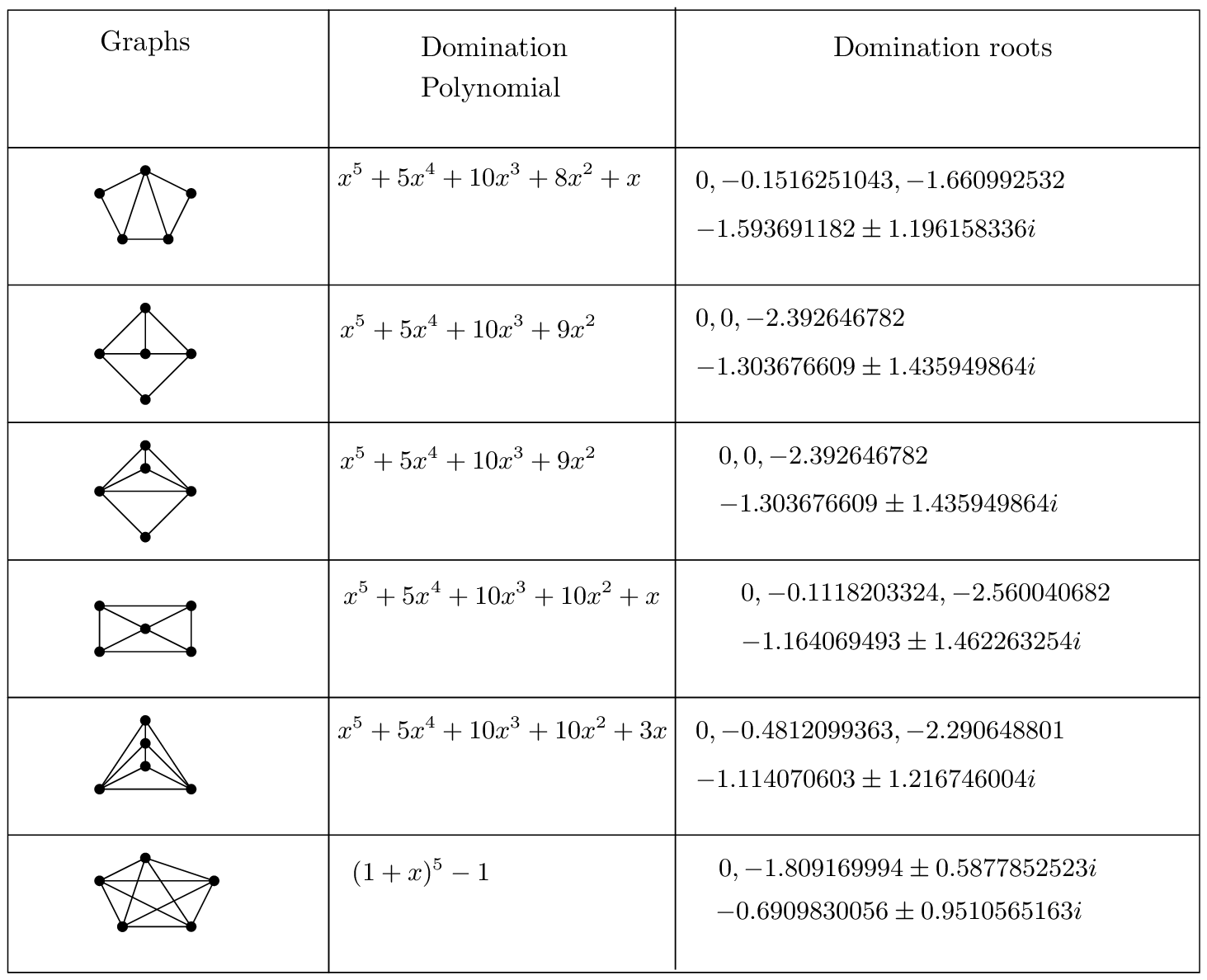}
\hglue0.75cm
\includegraphics[width=14cm,height=10.3cm]{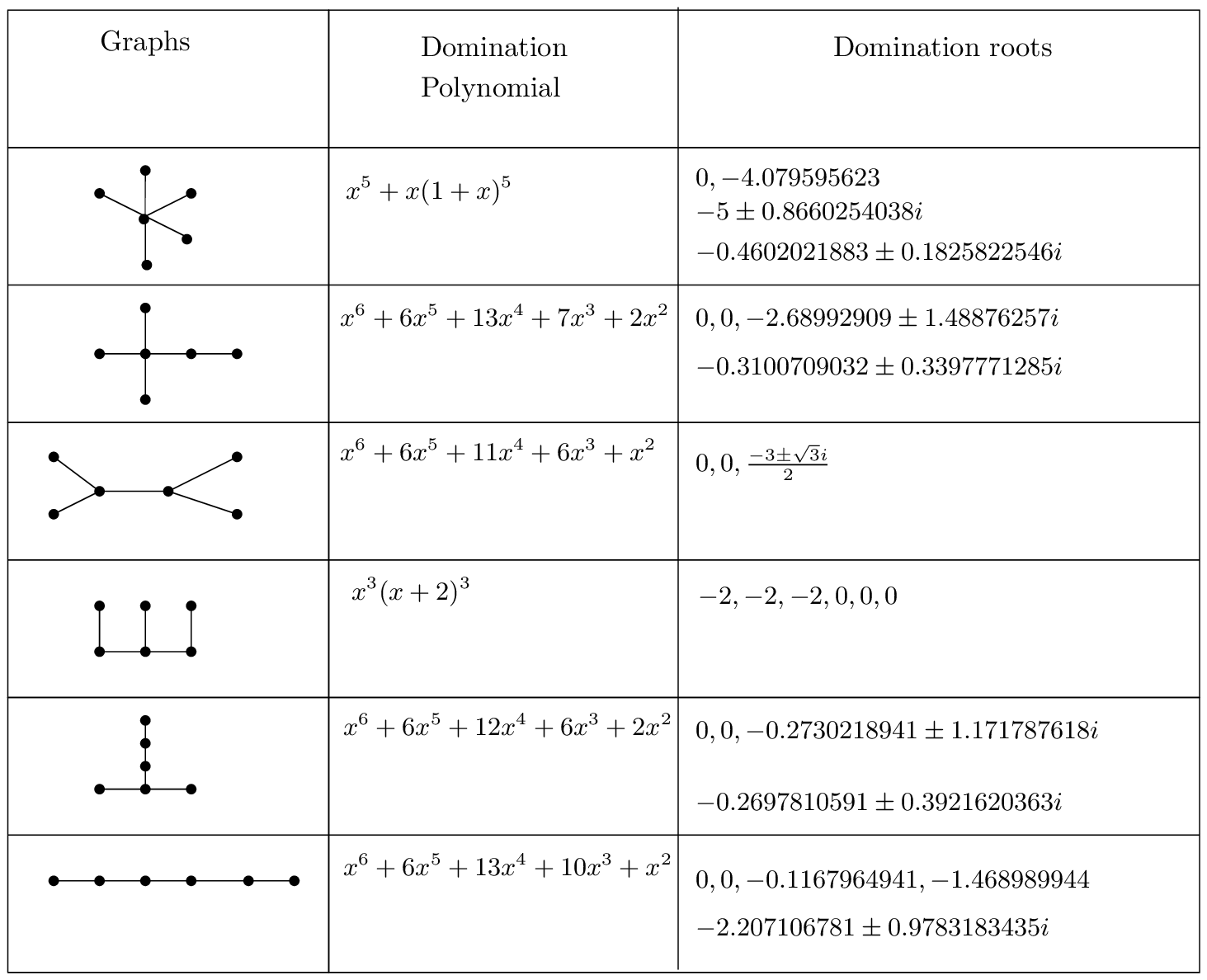}
\end{figure}

\begin{figure}[h]
\hglue.75cm
\includegraphics[width=14cm,height=10.3cm]{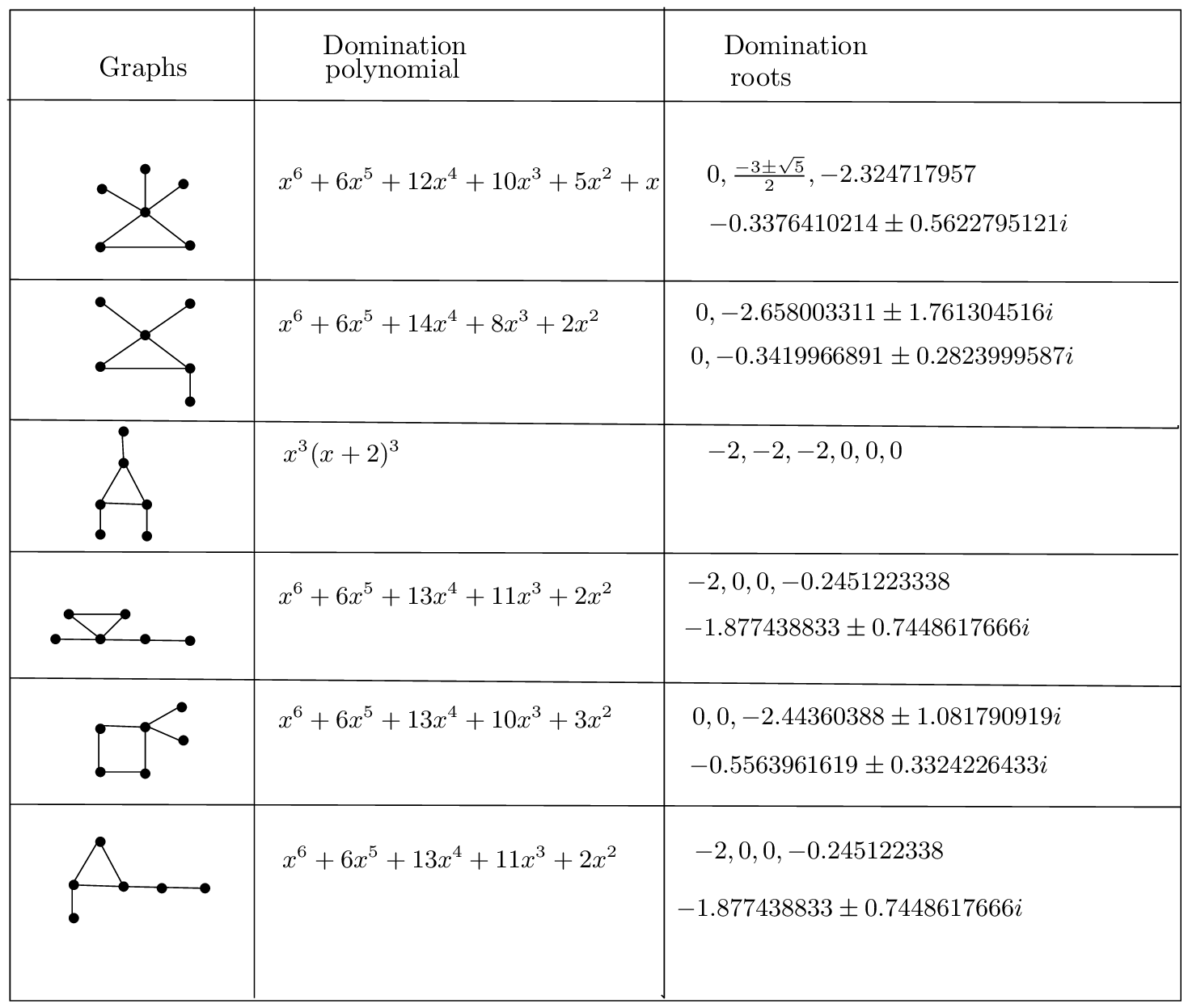}
\hglue.75cm
\includegraphics[width=14cm,height=10.3cm]{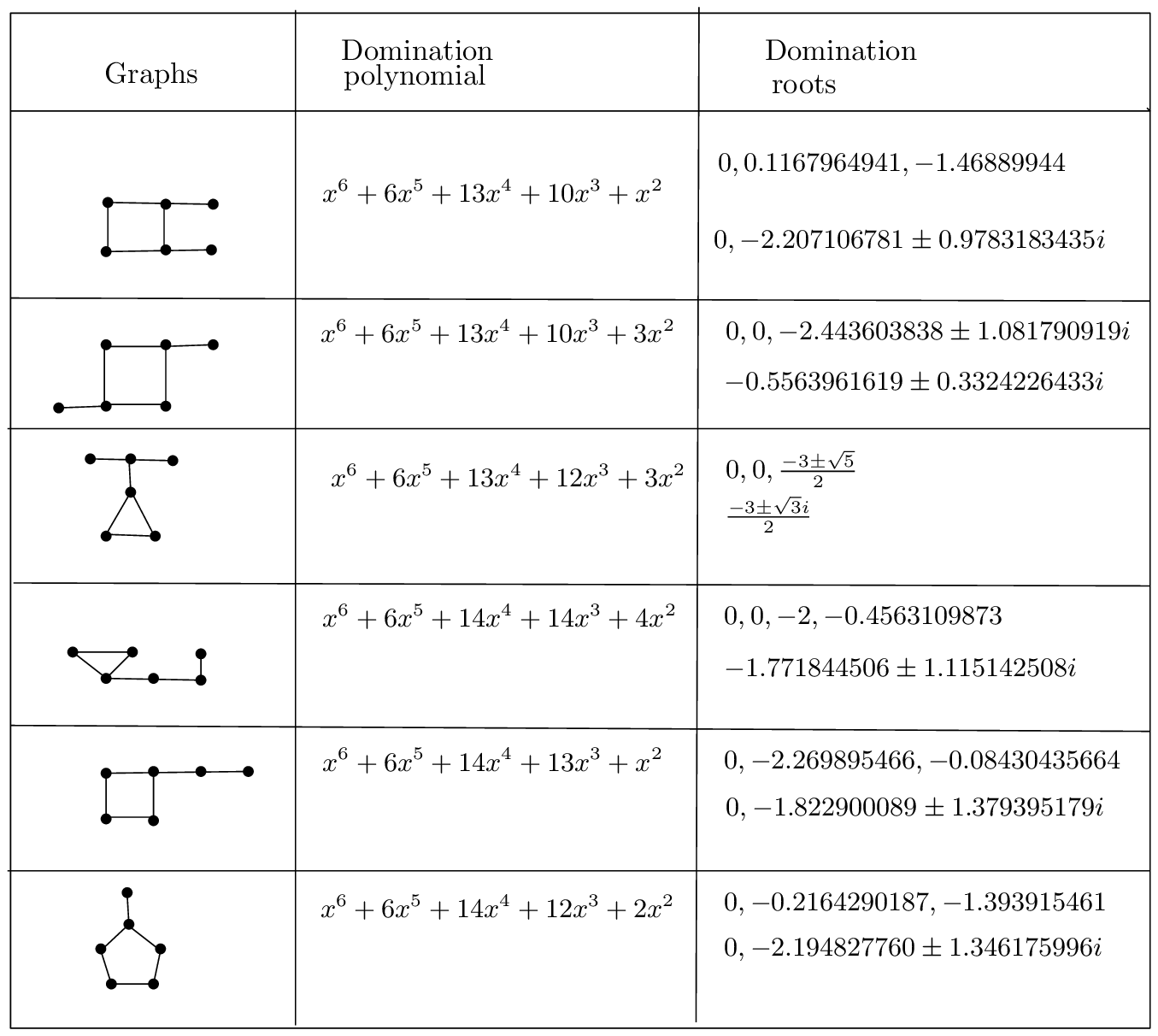}
\end{figure}

\begin{figure}[h]
\hglue0.75cm
\includegraphics[width=14cm,height=10.3cm]{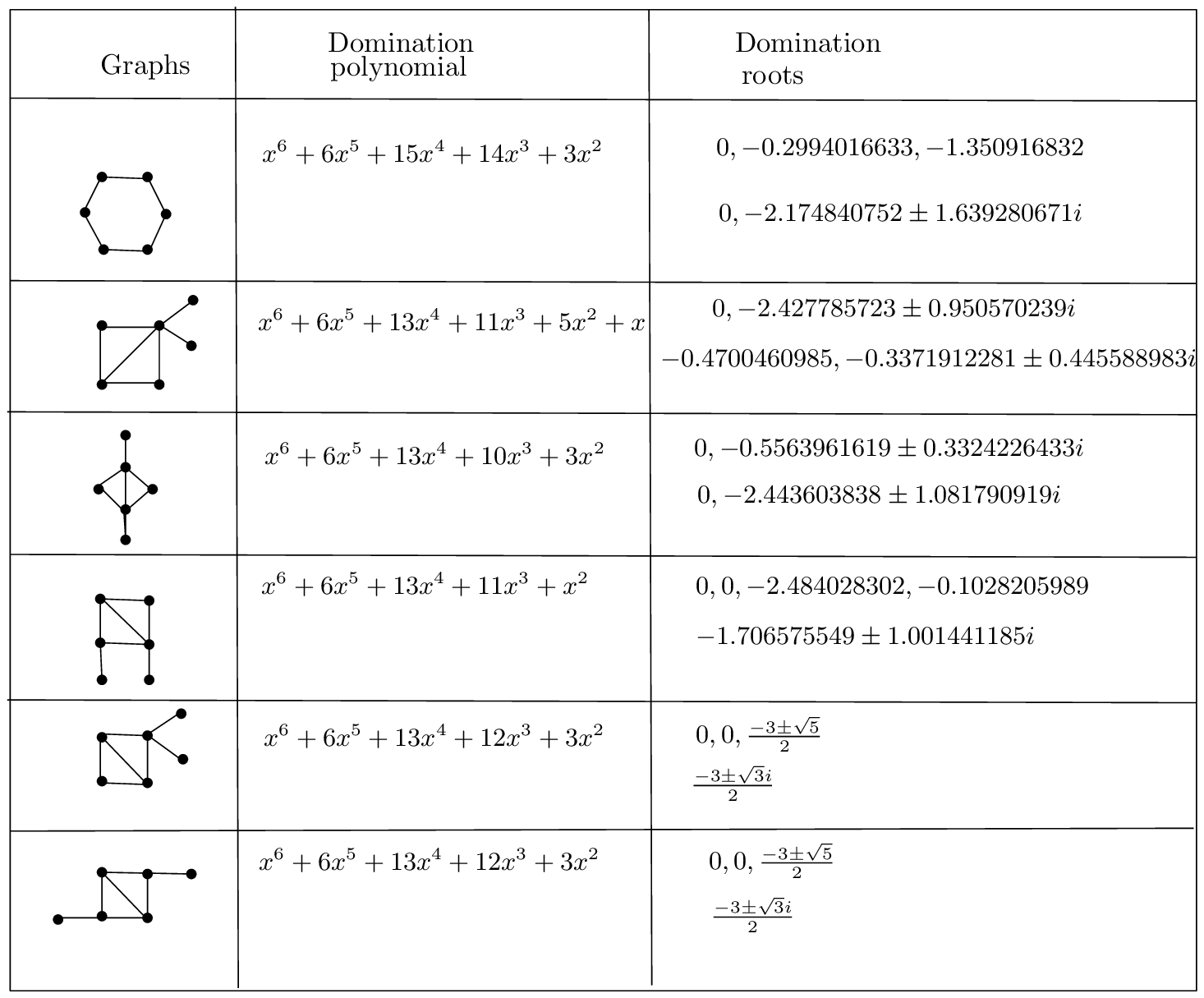}
\hglue0.75cm
\includegraphics[width=14cm,height=10.3cm]{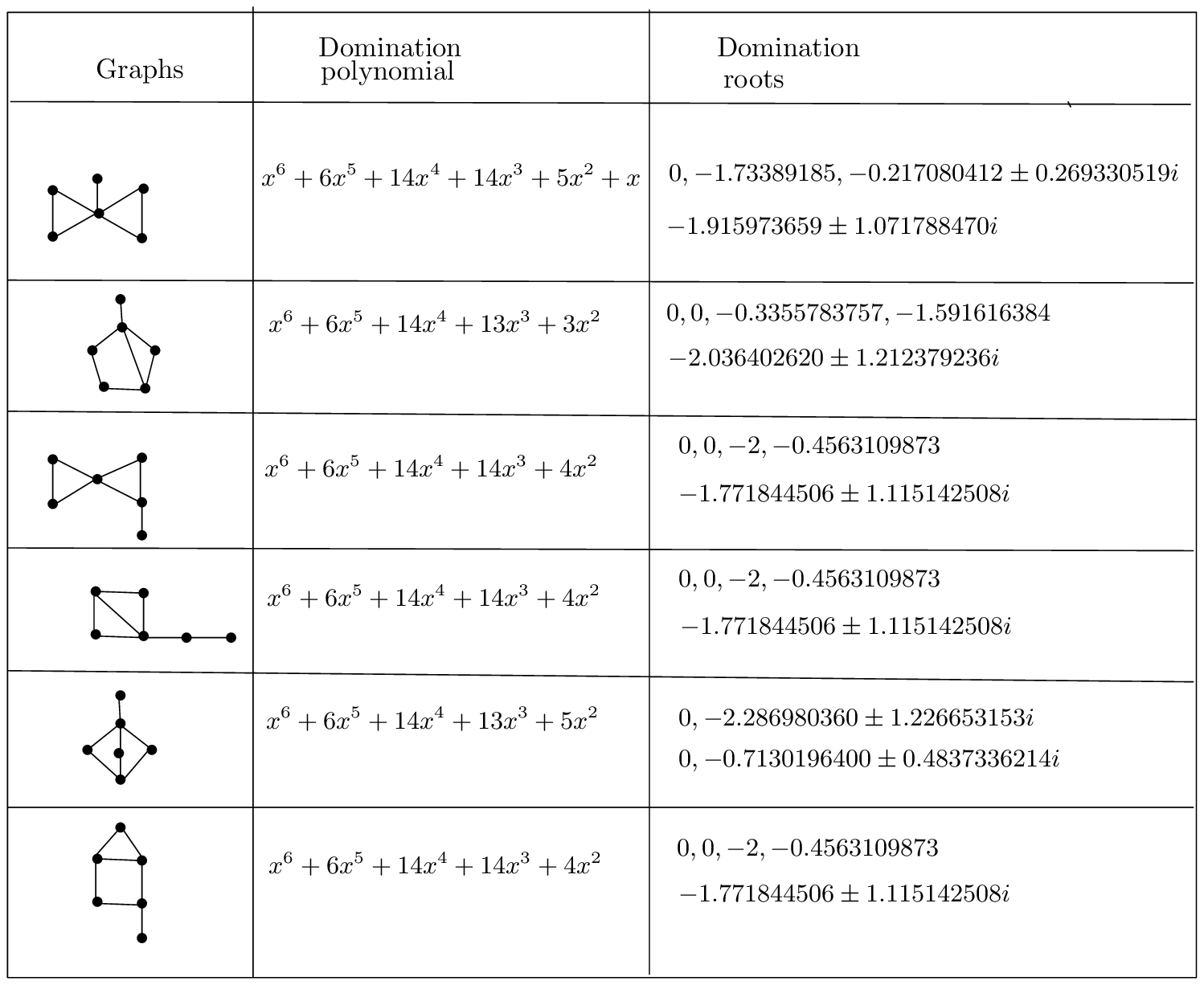}
\end{figure}

\begin{figure}[h]
\hglue0.75cm
\includegraphics[width=14cm,height=10.3cm]{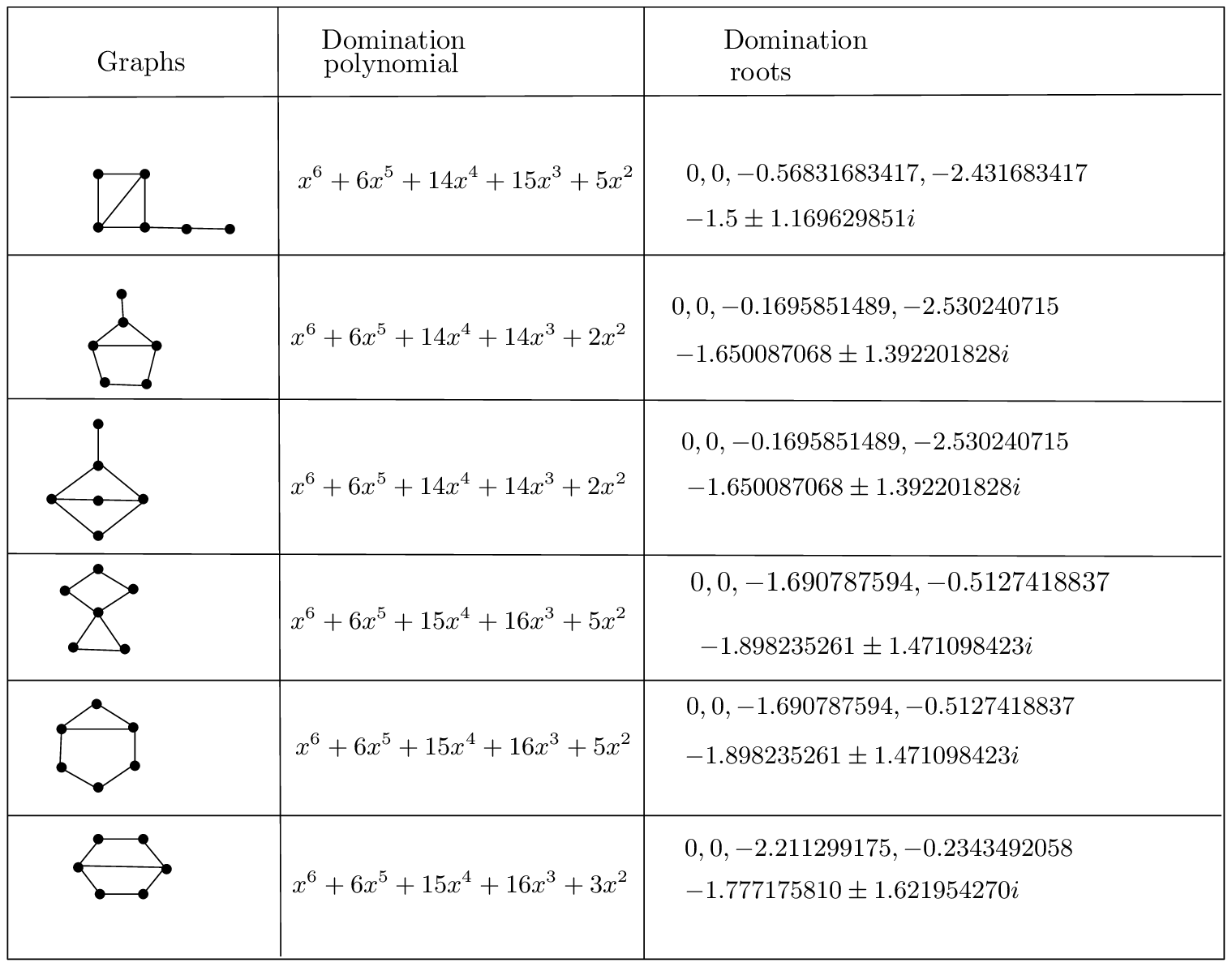}
\hglue0.55cm
\includegraphics[width=14.2cm,height=10.3cm]{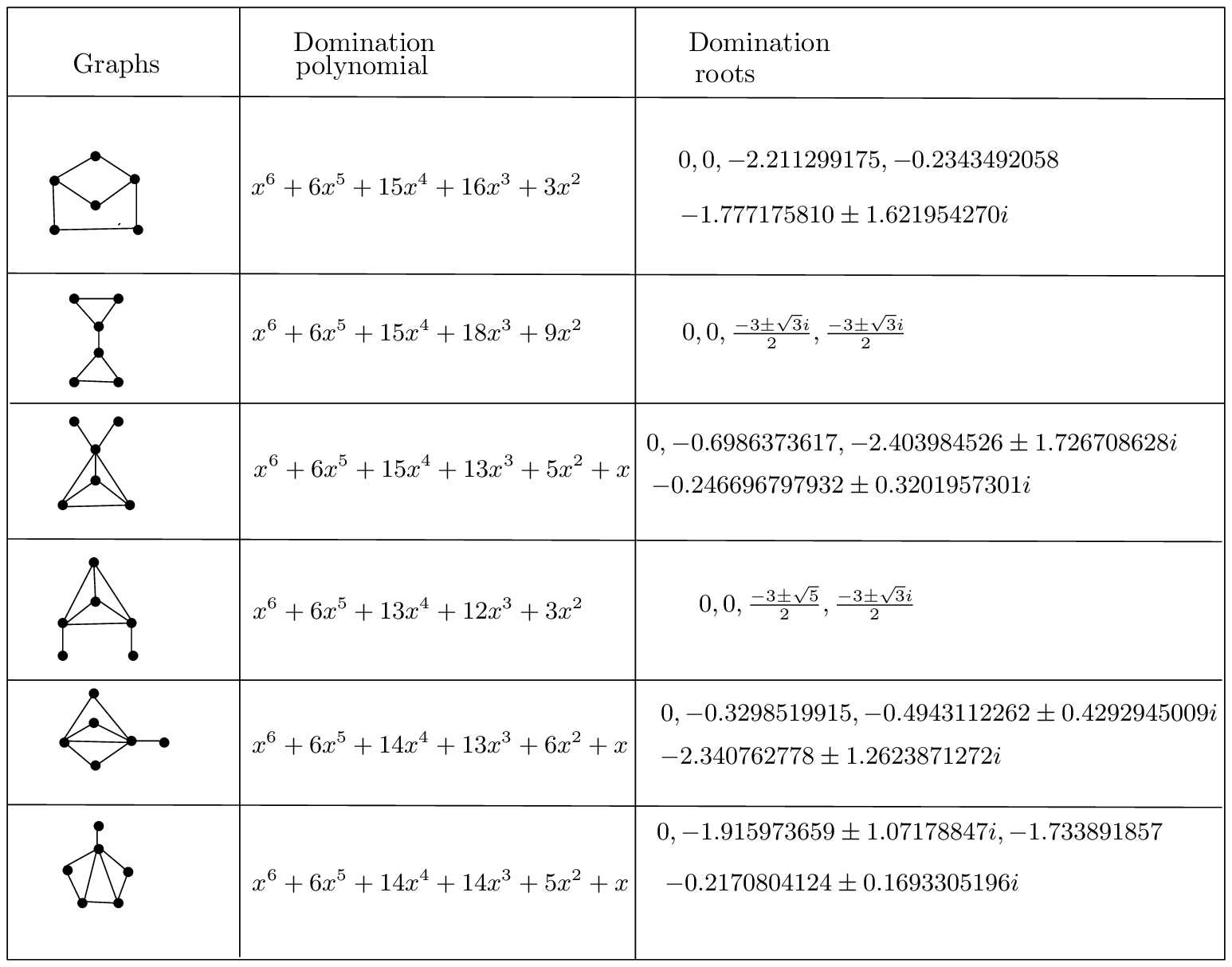}
\end{figure}

\begin{figure}[h]
\hglue0.60cm
\includegraphics[width=14.16cm,height=10.3cm]{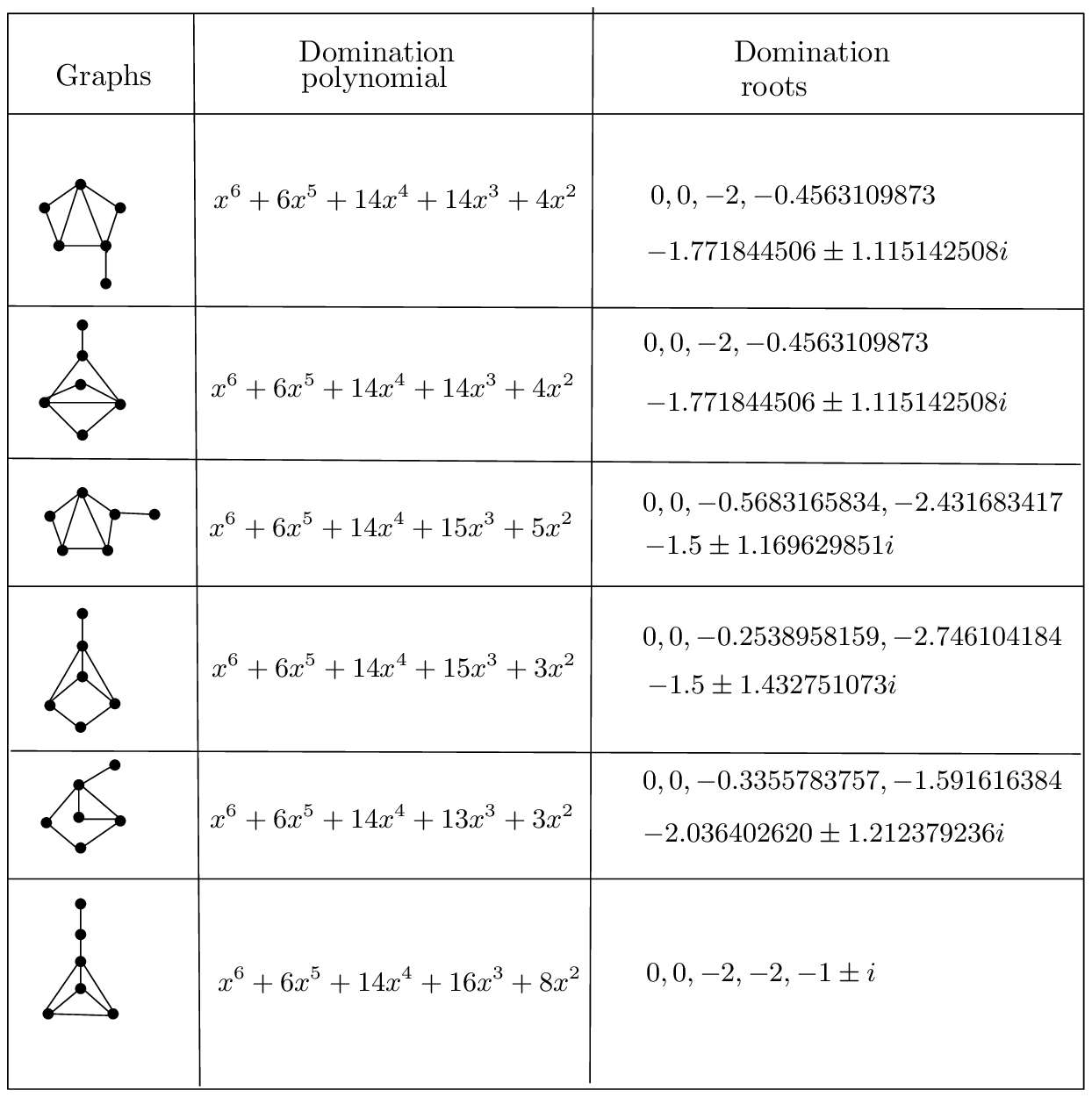}
\hglue.75cm
\includegraphics[width=14cm,height=10.3cm]{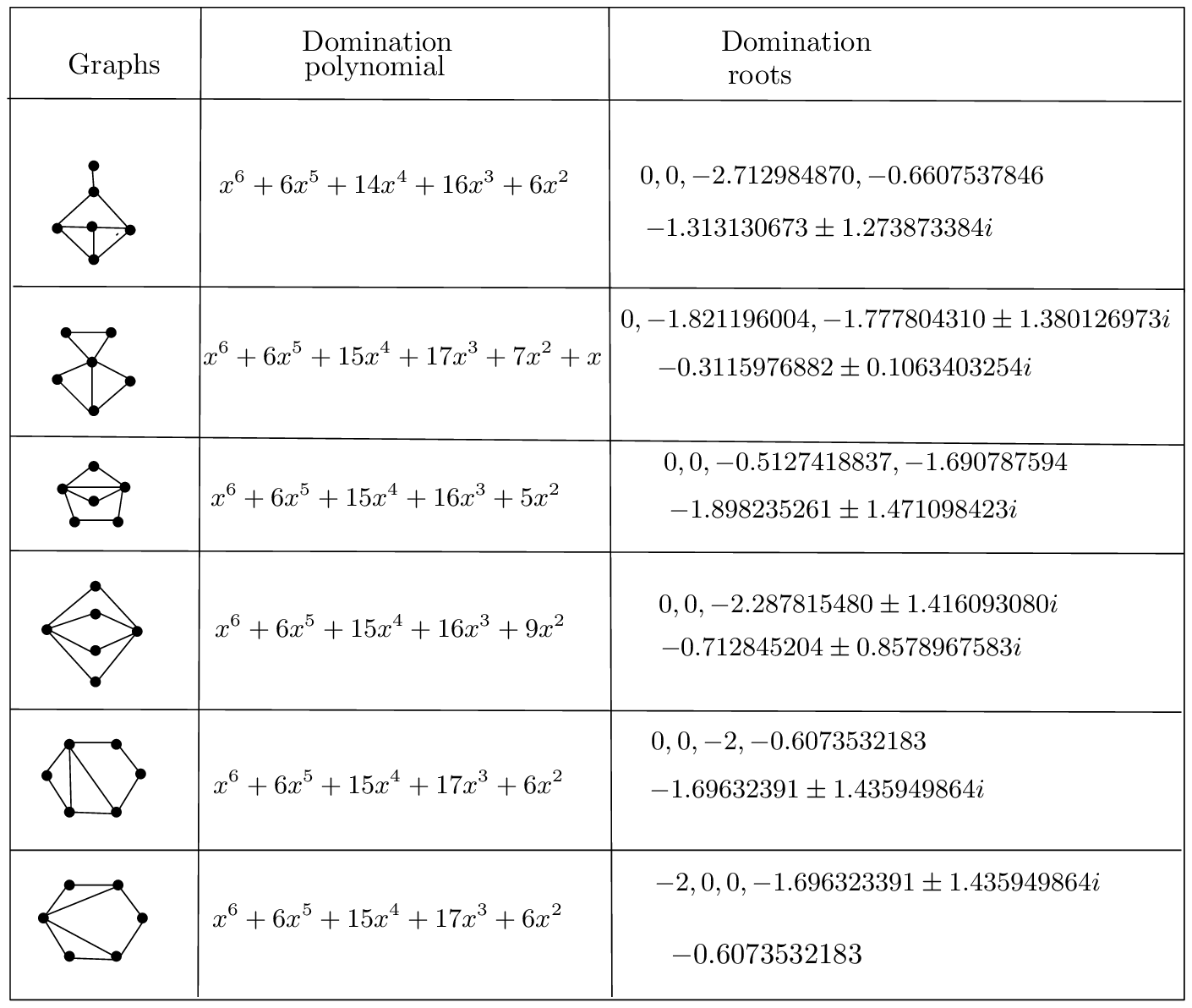}
\end{figure}

\begin{figure}
\hglue0.36cm
\includegraphics[width=14.4cm,height=10.3cm]{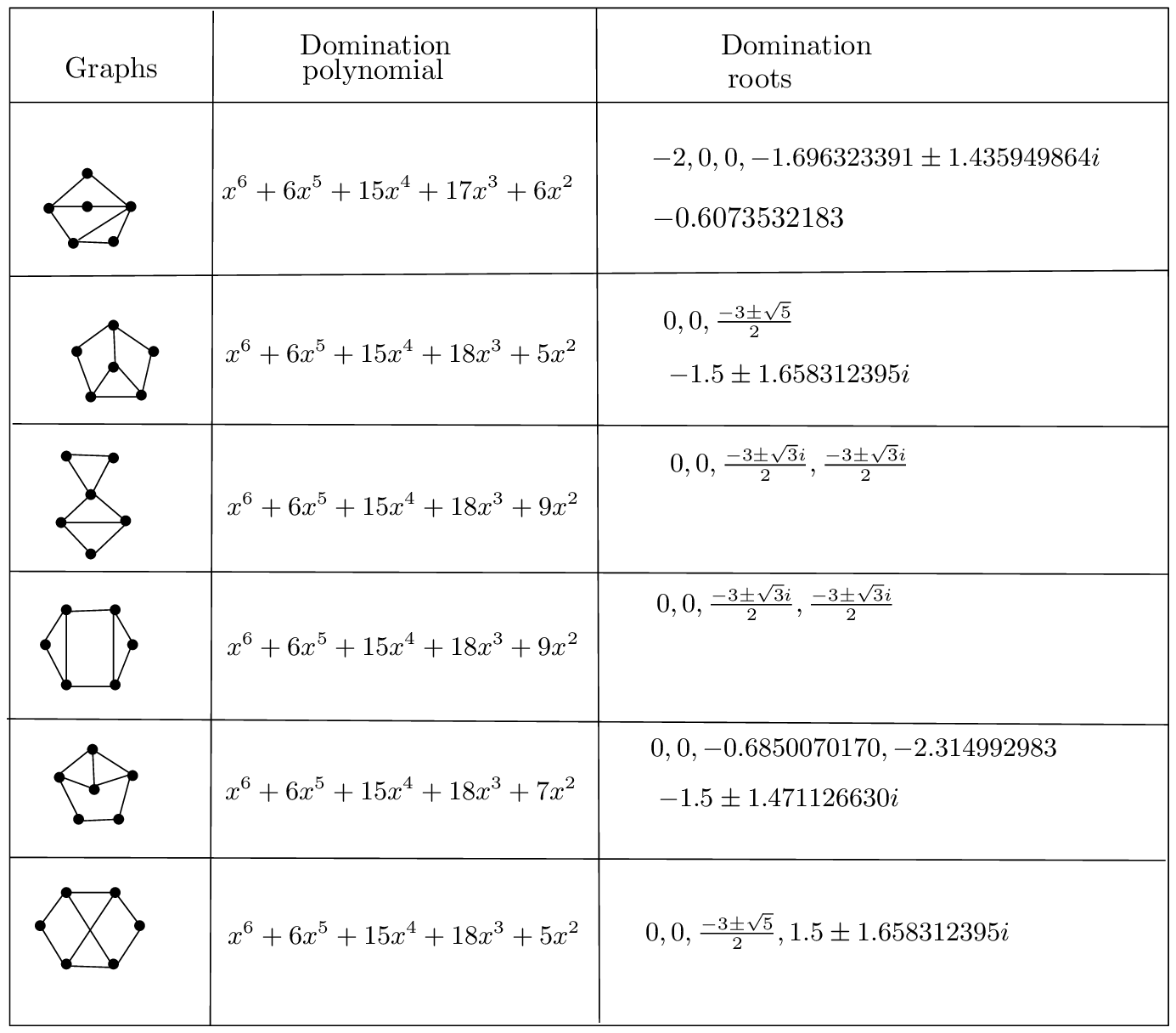}
\hglue0.75cm
\includegraphics[width=14cm,height=10.3cm]{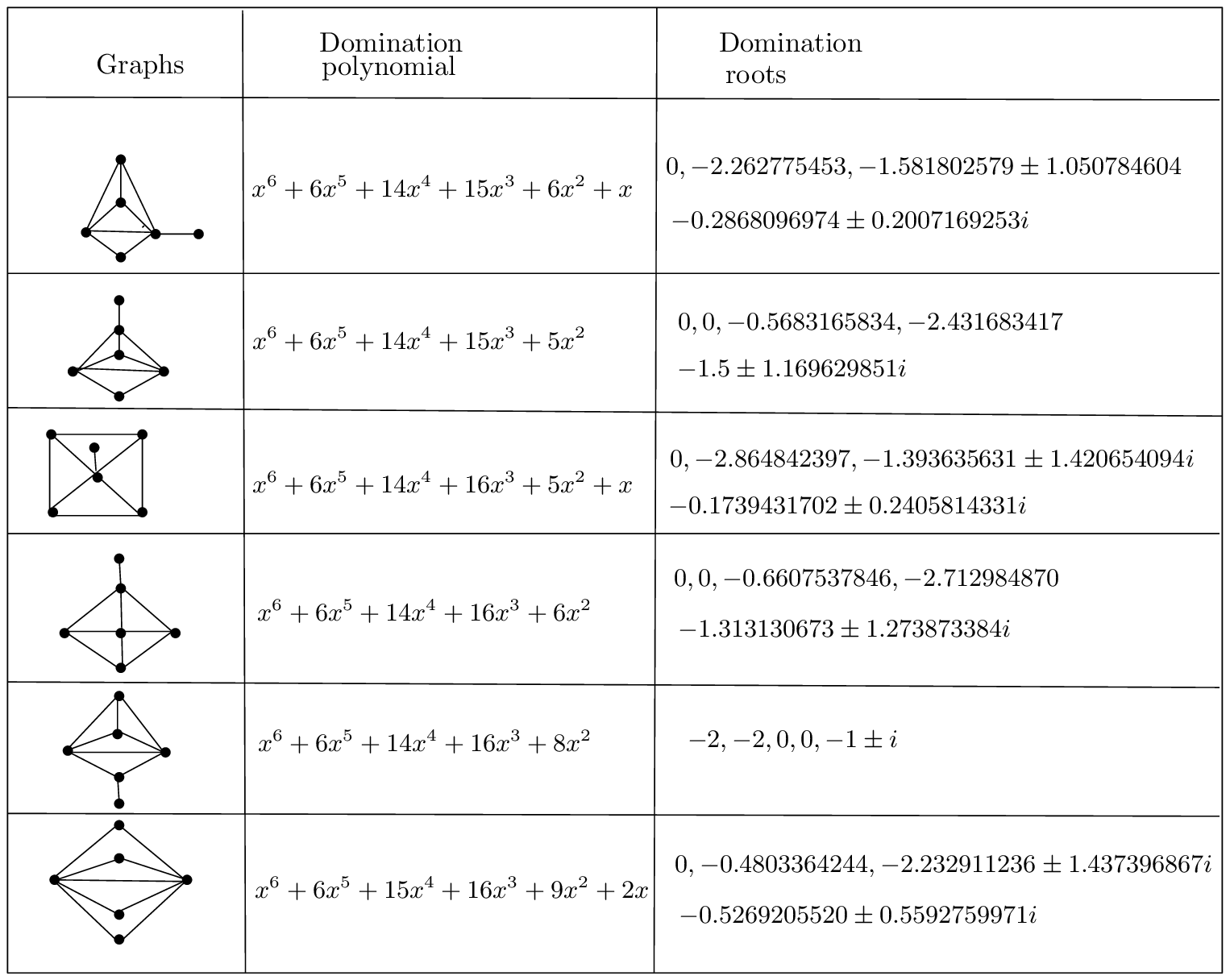}
\end{figure}
\begin{figure}[h]
\hglue0.6cm
\includegraphics[width=14.1cm,height=10.3cm]{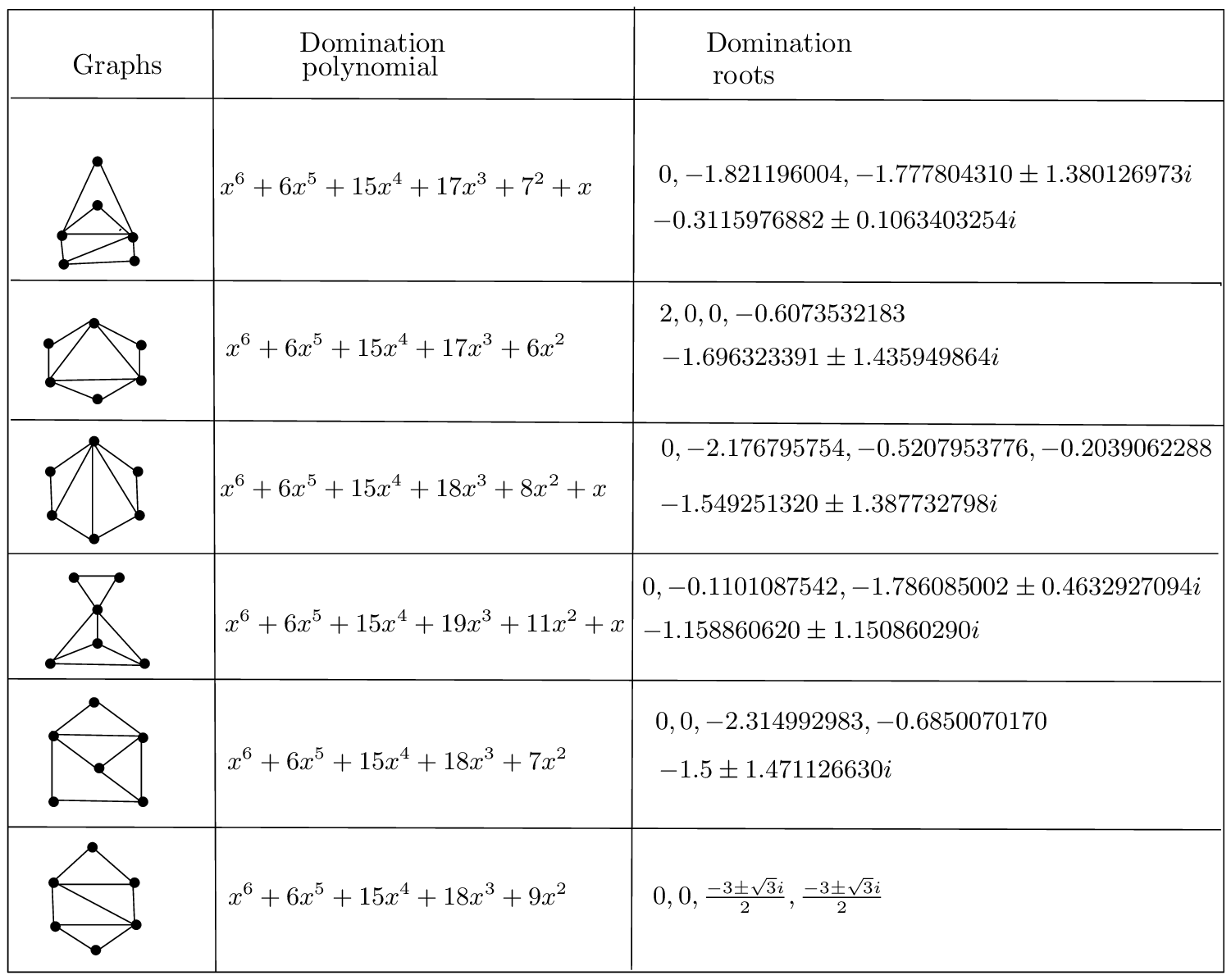}
\hglue0.8cm
\includegraphics[width=13.7cm,height=10.3cm]{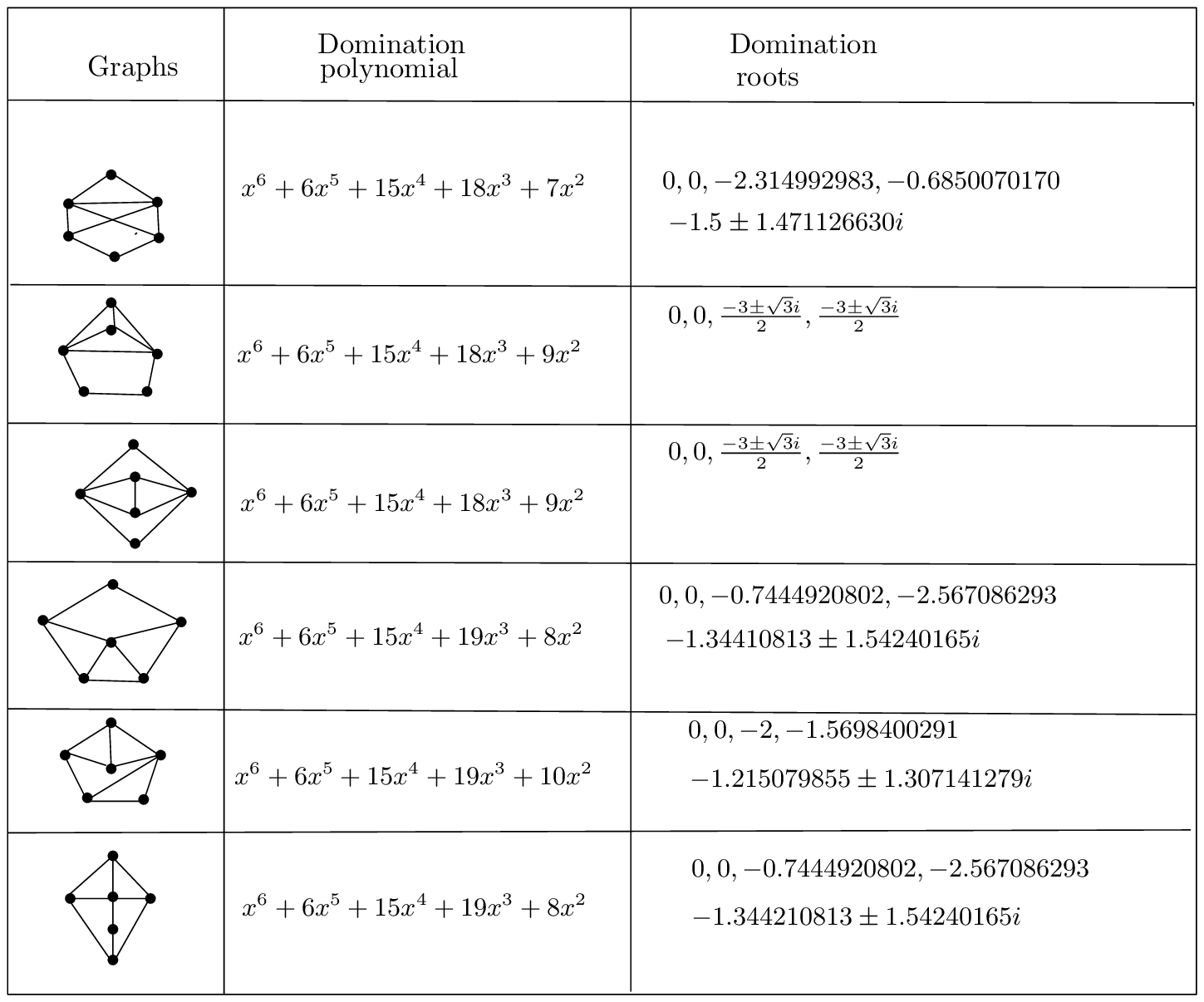}
\end{figure}
\begin{figure}
\hglue0.75cm
\includegraphics[width=14cm,height=10.3cm]{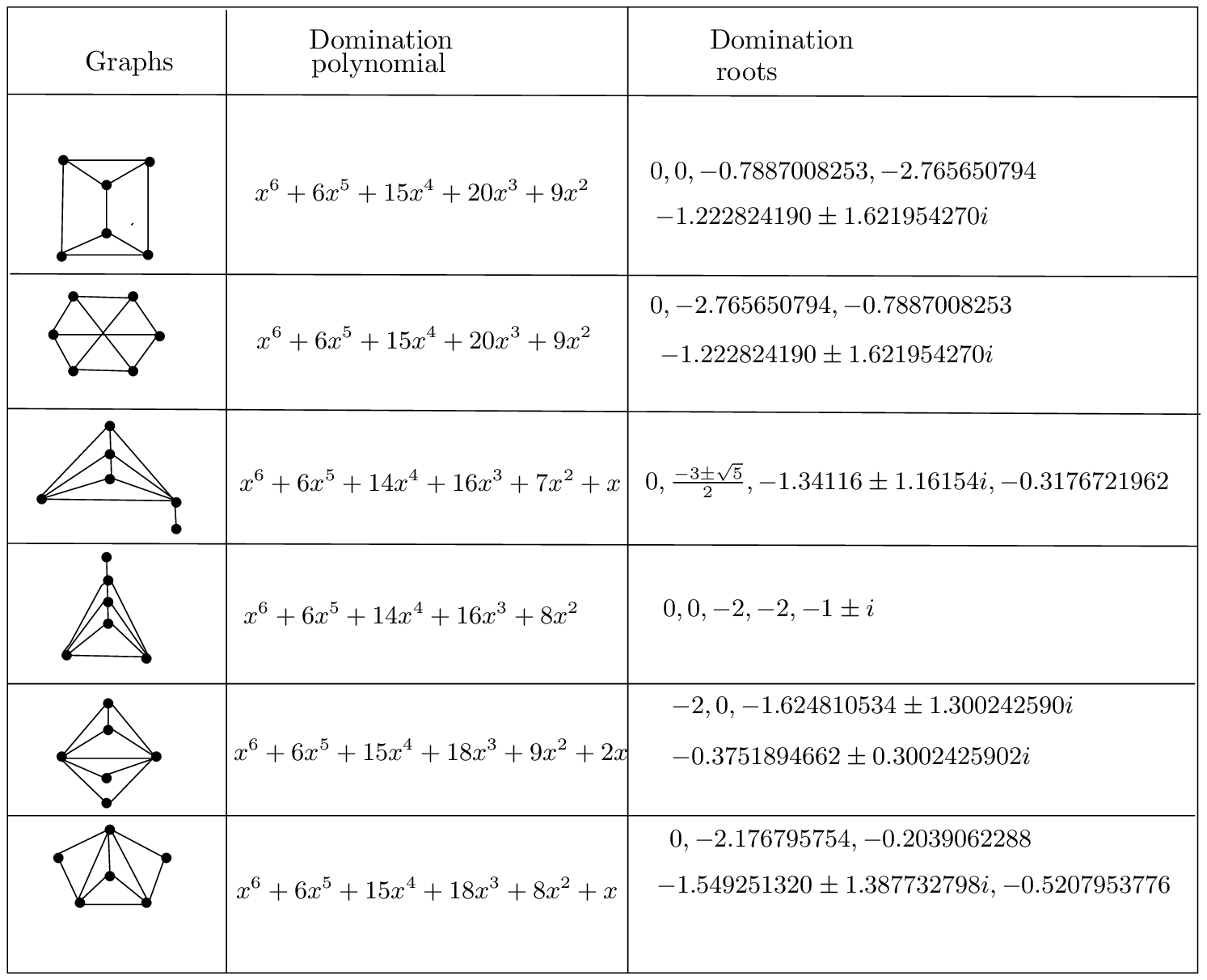}
\hglue0.75cm
\includegraphics[width=14cm,height=10.3cm]{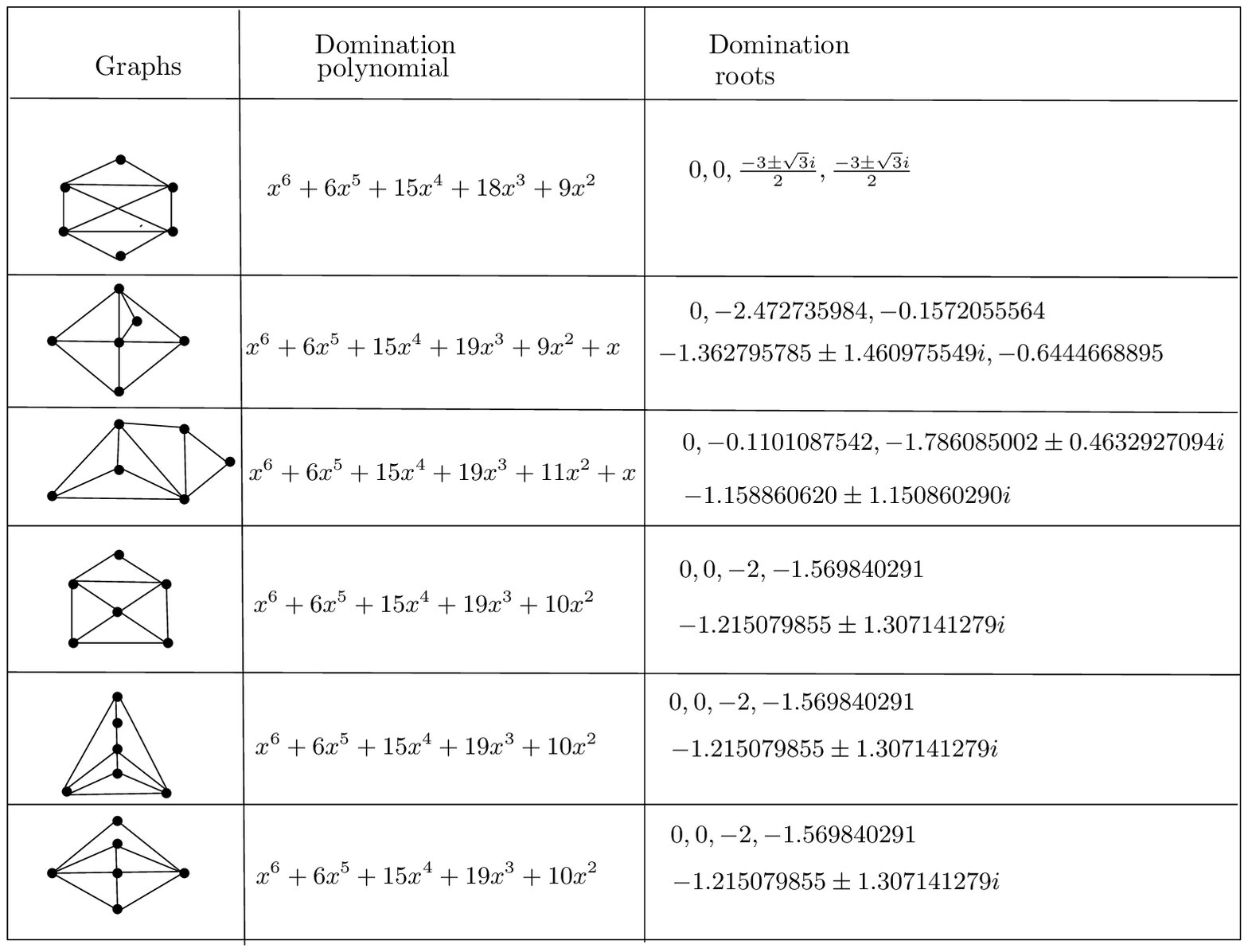}
\end{figure}

\begin{figure}[h]
\hglue0.75cm
\includegraphics[width=14cm,height=10.3cm]{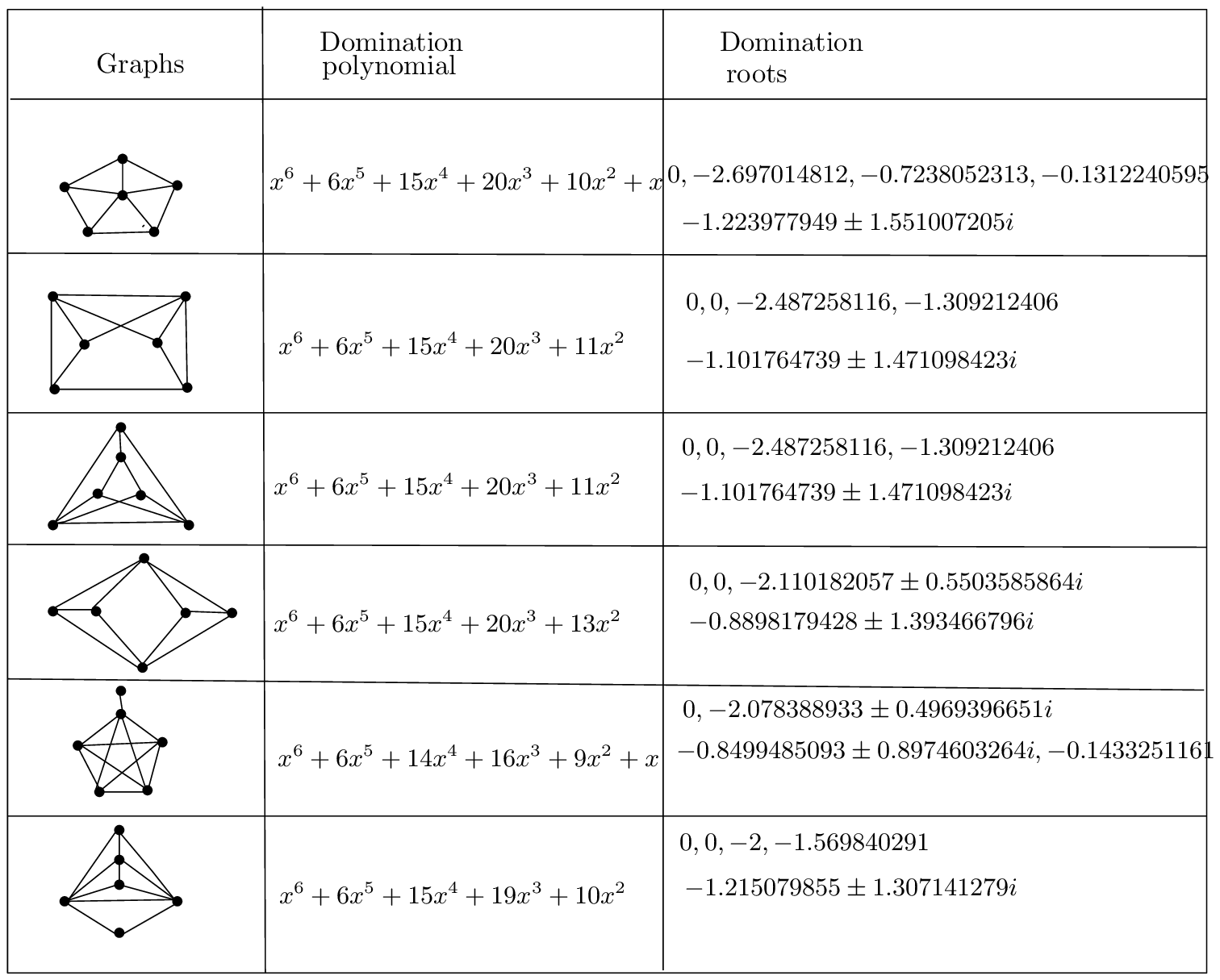}
\hglue0.75cm
\includegraphics[width=14cm,height=10.3cm]{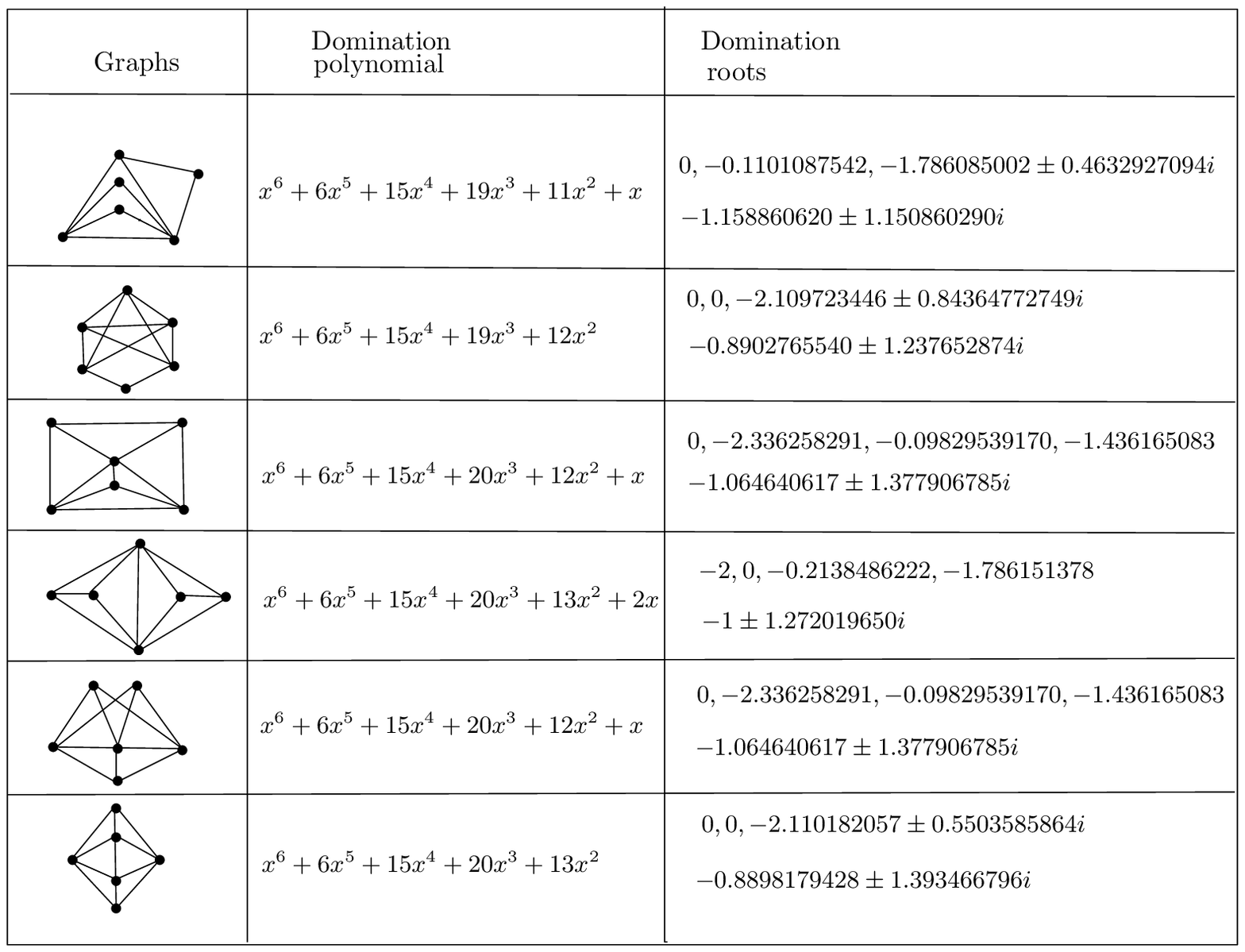}
\end{figure}
\begin{figure}
\hglue0.75cm
\includegraphics[width=14cm,height=10.3cm]{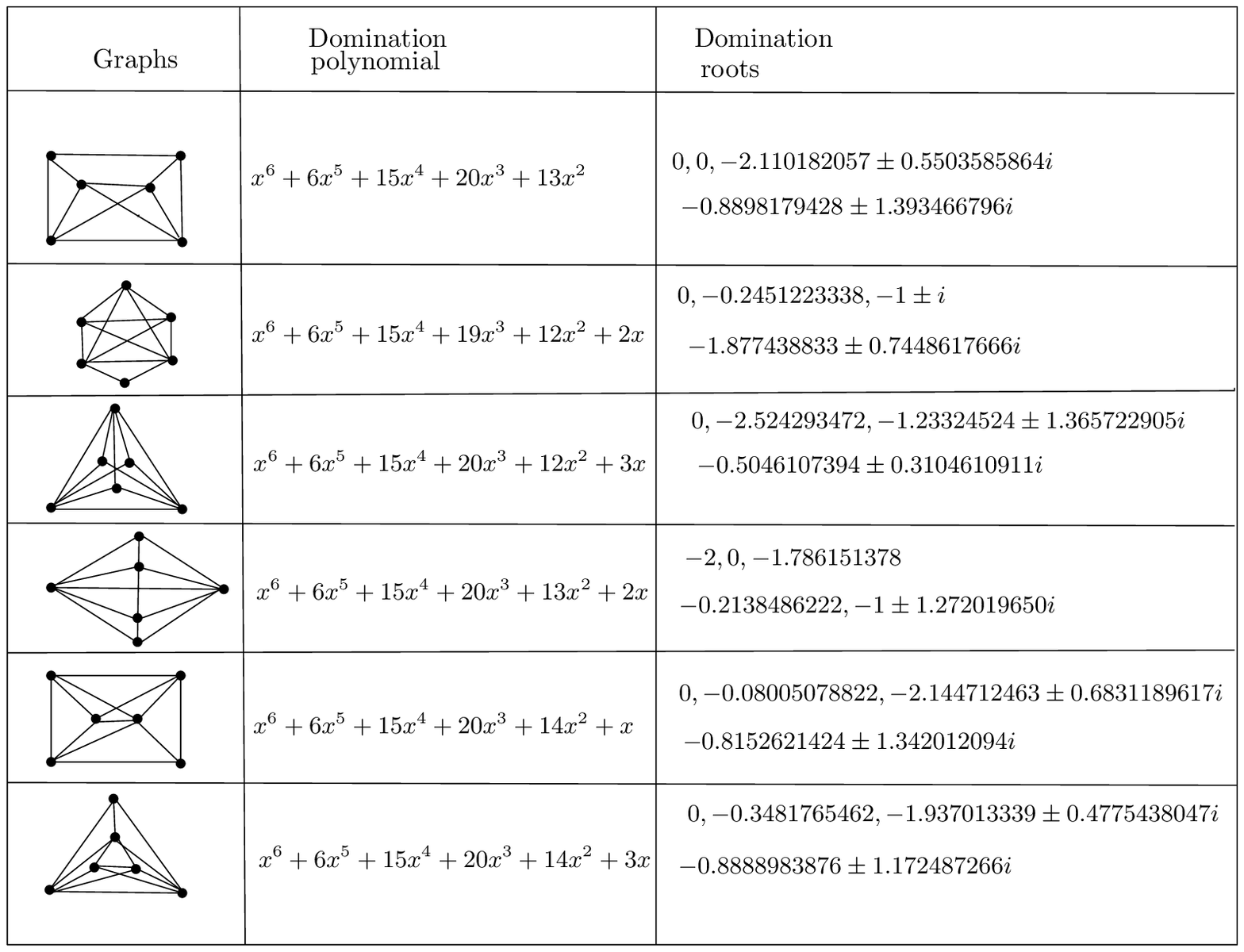}
\hglue0.75cm
\includegraphics[width=14cm,height=10.3cm]{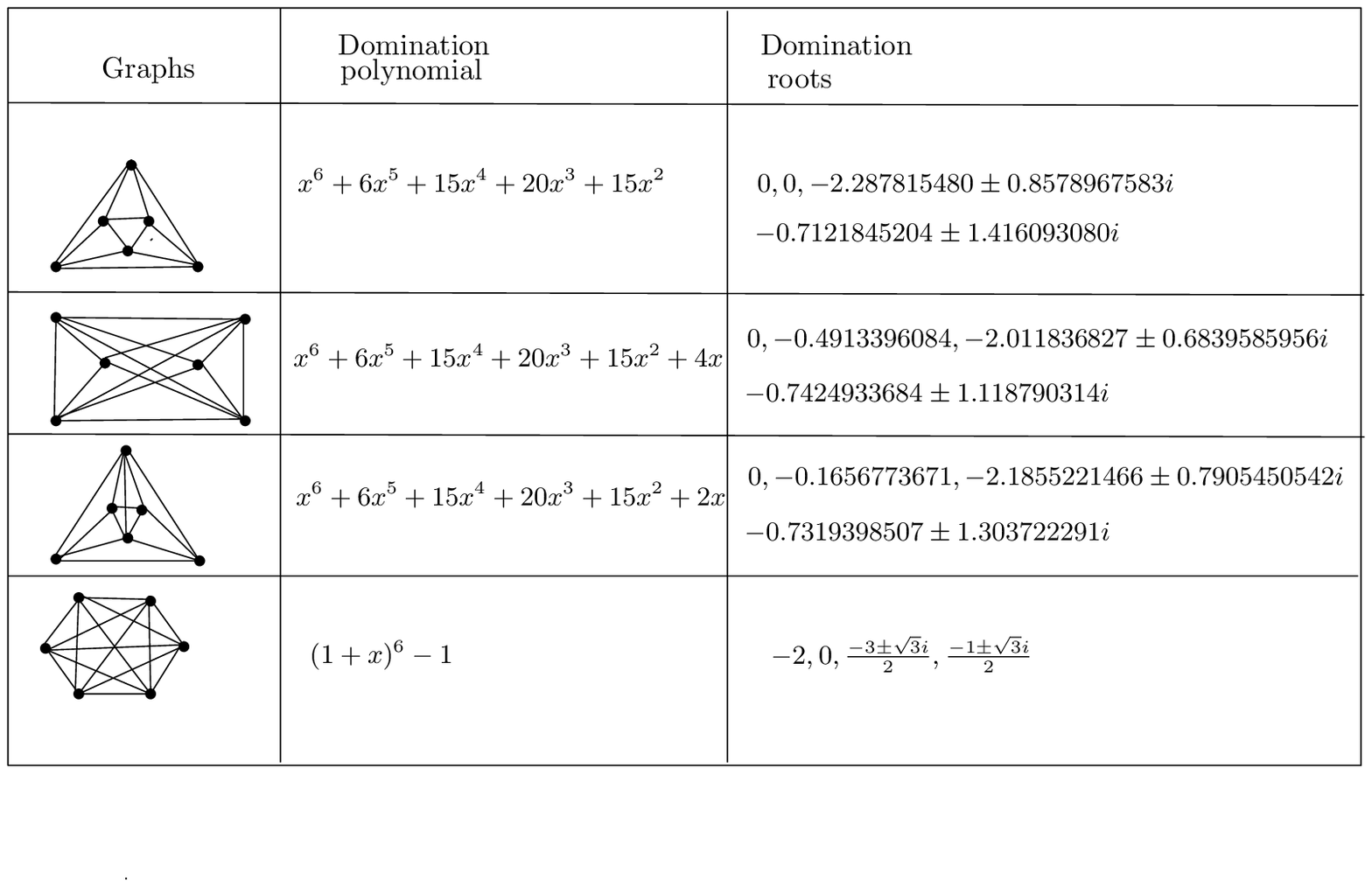}
\end{figure}




\begin{thebibliography}{99}




\bibitem{euro} S. Akbari, S. Alikhani and Y.H.  Peng, {\it Characterization of
graphs using domination polynomial}, Europ. J. Combin.,  Vol 31 (2010) 1714-1724.

\bibitem{Oper} S. Alikhani, {\it On the domination polynomial of some graph operations}, ISRN Combin.,
Volume 2013, Article ID 146595, 3 pages.

\bibitem{gcom} S. Alikhani, {\it The domination polynomial of a graph at $-1$}, Graphs  Combin., 29 (2013) 1175-1181.

\bibitem {saeid3}  S. Alikhani, {\it On the domination polynomials of non $P_4$-free graphs},
Iran. J. Math. Sci. Informatics, Vol. 8, No. 2 (2013) 49--55.


\bibitem{few} S. Alikhani, {\it Graphs whose certain polynomials have few distinct roots}, ISRN Discrete Math., Volume 2013,  Article ID 195818, 8 pages.

\bibitem{thesis} S. Alikhani, {\it Dominating sets and domination polynomials of graphs}, PhD thesis, Universiti Putra Malaysia, 2009. 

\bibitem{book} S. Alikhani, {\it Dominating sets and domination polynomials of graphs: Domination
polynomial: A new graph polynomial}, LAMBERT Academic Publishing, ISBN:
9783847344827 (2012).


\bibitem{saeid1} S. Alikhani, Y. H. Peng, {\it Introduction to domination polynomial of a graphs}, Ars Combin., to appear. Available at \texttt{http://arxiv.org/abs/0905.2251}.
\bibitem{domination} T.W. Haynes, S.T. Hedetniemi, P.J. Slater, {\it Fundamentals of Domination in Graphs, Marcel Dekker}, NewYork, 1998.
\bibitem{atlas} R.C. Read,  R.J. Wilson, {\it An atlas of graphs}, Clarendon Press, Oxford, 1998.
\end{thebibliography}
\end{document}